\documentclass[12pt]{amsart}
 \usepackage{amsmath}
\usepackage{amssymb,latexsym}
\usepackage{amscd}
\usepackage{comment}
 \usepackage{xspace}
 \usepackage{verbatim}
\usepackage[english]{babel}
\usepackage[pdftex,bookmarks=true]{hyperref}

\begin{document}

\numberwithin{equation}{section}
\title[Parabolic equations with exponential nonlinearity]{Parabolic equations with exponential nonlinearity and measure data}
\author{Phuoc-Tai Nguyen}
\address{Department of Mathematics, Technion\\
 Haifa 32000, ISRAEL}
 \email{nguyenphuoctai.hcmup@gmail.com}

\date{}

 \maketitle


\newcommand{\txt}[1]{\;\text{ #1 }\;}
\newcommand{\tbf}{\textbf}
\newcommand{\tit}{\textit}
\newcommand{\tsc}{\textsc}
\newcommand{\trm}{\textrm}
\newcommand{\mbf}{\mathbf}
\newcommand{\mrm}{\mathrm}
\newcommand{\bsym}{\boldsymbol}
\newcommand{\scs}{\scriptstyle}
\newcommand{\sss}{\scriptscriptstyle}
\newcommand{\txts}{\textstyle}
\newcommand{\dsps}{\displaystyle}
\newcommand{\fnz}{\footnotesize}
\newcommand{\scz}{\scriptsize}
\newcommand{\be}{\begin{equation}}
\newcommand{\bel}[1]{\begin{equation}\label{#1}}
\newcommand{\ee}{\end{equation}}
\newcommand{\eqnl}[2]{\begin{equation}\label{#1}{#2}\end{equation}}
\newcommand{\barr}{\begin{eqnarray}}
\newcommand{\earr}{\end{eqnarray}}
\newcommand{\bars}{\begin{eqnarray*}}
\newcommand{\ears}{\end{eqnarray*}}
\newcommand{\nnu}{\nonumber \\}
\newtheorem{subn}{\name}
\renewcommand{\thesubn}{}
\newcommand{\bsn}[1]{\def\name{#1}\begin{subn}}
\newcommand{\esn}{\end{subn}}
\newtheorem{sub}{\name}[section]
\newcommand{\dn}[1]{\def\name{#1}}   
\newcommand{\bs}{\begin{sub}}
\newcommand{\es}{\end{sub}}
\newcommand{\bsl}[1]{\begin{sub}\label{#1}}
\newcommand{\bth}[1]{\def\name{Theorem}
\begin{sub}\label{t:#1}}
\newcommand{\blemma}[1]{\def\name{Lemma}
\begin{sub}\label{l:#1}}
\newcommand{\bcor}[1]{\def\name{Corollary}
\begin{sub}\label{c:#1}}
\newcommand{\bdef}[1]{\def\name{Definition}
\begin{sub}\label{d:#1}}
\newcommand{\bprop}[1]{\def\name{Proposition}
\begin{sub}\label{p:#1}}
\newcommand{\R}{\eqref}
\newcommand{\rth}[1]{Theorem~\ref{t:#1}}
\newcommand{\rlemma}[1]{Lemma~\ref{l:#1}}
\newcommand{\rcor}[1]{Corollary~\ref{c:#1}}
\newcommand{\rdef}[1]{Definition~\ref{d:#1}}
\newcommand{\rprop}[1]{Proposition~\ref{p:#1}}
\newcommand{\BA}{\begin{array}}
\newcommand{\EA}{\end{array}}
\newcommand{\BAN}{\renewcommand{\arraystretch}{1.2}
\setlength{\arraycolsep}{2pt}\begin{array}}
\newcommand{\BAV}[2]{\renewcommand{\arraystretch}{#1}
\setlength{\arraycolsep}{#2}\begin{array}}
\newcommand{\BSA}{\begin{subarray}}
\newcommand{\ESA}{\end{subarray}}
\newcommand{\BAL}{\begin{aligned}}
\newcommand{\EAL}{\end{aligned}}
\newcommand{\BALG}{\begin{alignat}}
\newcommand{\EALG}{\end{alignat}}
\newcommand{\BALGN}{\begin{alignat*}}
\newcommand{\EALGN}{\end{alignat*}}
\newcommand{\note}[1]{\textit{#1.}\hspace{2mm}}
\newcommand{\Proof}{\note{Proof}}
\newcommand{\Remark}{\note{Remark}}
\newcommand{\modin}{$\,$\\[-4mm] \indent}
\newcommand{\forevery}{\quad \forall}
\newcommand{\set}[1]{\{#1\}}
\newcommand{\setdef}[2]{\{\,#1:\,#2\,\}}
\newcommand{\setm}[2]{\{\,#1\mid #2\,\}}
\newcommand{\mt}{\mapsto}
\newcommand{\lra}{\longrightarrow}
\newcommand{\lla}{\longleftarrow}
\newcommand{\llra}{\longleftrightarrow}
\newcommand{\Lra}{\Longrightarrow}
\newcommand{\Lla}{\Longleftarrow}
\newcommand{\Llra}{\Longleftrightarrow}
\newcommand{\warrow}{\rightharpoonup}
\newcommand{
\paran}[1]{\left (#1 \right )}
\newcommand{\sqbr}[1]{\left [#1 \right ]}
\newcommand{\curlybr}[1]{\left \{#1 \right \}}
\newcommand{\abs}[1]{\left |#1\right |}
\newcommand{\norm}[1]{\left \|#1\right \|}
\newcommand{
\paranb}[1]{\big (#1 \big )}
\newcommand{\lsqbrb}[1]{\big [#1 \big ]}
\newcommand{\lcurlybrb}[1]{\big \{#1 \big \}}
\newcommand{\absb}[1]{\big |#1\big |}
\newcommand{\normb}[1]{\big \|#1\big \|}
\newcommand{
\paranB}[1]{\Big (#1 \Big )}
\newcommand{\absB}[1]{\Big |#1\Big |}
\newcommand{\normB}[1]{\Big \|#1\Big \|}
\newcommand{\produal}[1]{\langle #1 \rangle}

\newcommand{\thkl}{\rule[-.5mm]{.3mm}{3mm}}
\newcommand{\thknorm}[1]{\thkl #1 \thkl\,}
\newcommand{\trinorm}[1]{|\!|\!| #1 |\!|\!|\,}
\newcommand{\bang}[1]{\langle #1 \rangle}
\def\angb<#1>{\langle #1 \rangle}
\newcommand{\vstrut}[1]{\rule{0mm}{#1}}
\newcommand{\rec}[1]{\frac{1}{#1}}
\newcommand{\opname}[1]{\mbox{\rm #1}\,}
\newcommand{\supp}{\opname{supp}}
\newcommand{\dist}{\opname{dist}}
\newcommand{\myfrac}[2]{{\displaystyle \frac{#1}{#2} }}
\newcommand{\myint}[2]{{\displaystyle \int_{#1}^{#2}}}
\newcommand{\mysum}[2]{{\displaystyle \sum_{#1}^{#2}}}
\newcommand {\dint}{{\displaystyle \myint\!\!\myint}}
\newcommand{\q}{\quad}
\newcommand{\qq}{\qquad}
\newcommand{\hsp}[1]{\hspace{#1mm}}
\newcommand{\vsp}[1]{\vspace{#1mm}}
\newcommand{\ity}{\infty}
\newcommand{\prt}{\partial}
\newcommand{\sms}{\setminus}
\newcommand{\ems}{\emptyset}
\newcommand{\ti}{\times}
\newcommand{\pr}{^\prime}
\newcommand{\ppr}{^{\prime\prime}}
\newcommand{\tl}{\tilde}
\newcommand{\sbs}{\subset}
\newcommand{\sbeq}{\subseteq}
\newcommand{\nind}{\noindent}
\newcommand{\ind}{\indent}
\newcommand{\ovl}{\overline}
\newcommand{\unl}{\underline}
\newcommand{\nin}{\not\in}
\newcommand{\pfrac}[2]{\genfrac{(}{)}{}{}{#1}{#2}}

\def\ga{\alpha}     \def\gb{\beta}       \def\gg{\gamma}
\def\gc{\chi}       \def\gd{\delta}      \def\ge{\epsilon}
\def\gth{\theta}                         \def\vge{\varepsilon}
\def\gf{\phi}       \def\vgf{\varphi}    \def\gh{\eta}
\def\gi{\iota}      \def\gk{\kappa}      \def\gl{\lambda}
\def\gm{\mu}        \def\gn{\nu}         \def\gp{\pi}
\def\vgp{\varpi}    \def\gr{\rho}        \def\vgr{\varrho}
\def\gs{\sigma}     \def\vgs{\varsigma}  \def\gt{\tau}
\def\gu{\upsilon}   \def\gv{\vartheta}   \def\gw{\omega}
\def\gx{\xi}        \def\gy{\psi}        \def\gz{\zeta}
\def\Gg{\Gamma}     \def\Gd{\Delta}      \def\Gf{\Phi}
\def\Gth{\Theta}
\def\Gl{\Lambda}    \def\Gs{\Sigma}      \def\Gp{\Pi}
\def\Gw{\Omega}     \def\Gx{\Xi}         \def\Gy{\Psi}

\def\CS{{\mathcal S}}   \def\CM{{\mathcal M}}   \def\CN{{\mathcal N}}
\def\CR{{\mathcal R}}   \def\CO{{\mathcal O}}   \def\CP{{\mathcal P}}
\def\CA{{\mathcal A}}   \def\CB{{\mathcal B}}   \def\CC{{\mathcal C}}
\def\CD{{\mathcal D}}   \def\CE{{\mathcal E}}   \def\CF{{\mathcal F}}
\def\CG{{\mathcal G}}   \def\CH{{\mathcal H}}   \def\CI{{\mathcal I}}
\def\CJ{{\mathcal J}}   \def\CK{{\mathcal K}}   \def\CL{{\mathcal L}}
\def\CT{{\mathcal T}}   \def\CU{{\mathcal U}}   \def\CV{{\mathcal V}}
\def\CZ{{\mathcal Z}}   \def\CX{{\mathcal X}}   \def\CY{{\mathcal Y}}
\def\CW{{\mathcal W}} \def\CQ{{\mathcal Q}}
\def\BBA {\mathbb A}   \def\BBb {\mathbb B}    \def\BBC {\mathbb C}
\def\BBD {\mathbb D}   \def\BBE {\mathbb E}    \def\BBF {\mathbb F}
\def\BBG {\mathbb G}   \def\BBH {\mathbb H}    \def\BBI {\mathbb I}
\def\BBJ {\mathbb J}   \def\BBK {\mathbb K}    \def\BBL {\mathbb L}
\def\BBM {\mathbb M}   \def\BBN {\mathbb N}    \def\BBO {\mathbb O}
\def\BBP {\mathbb P}   \def\BBR {\mathbb R}    \def\BBS {\mathbb S}
\def\BBT {\mathbb T}   \def\BBU {\mathbb U}    \def\BBV {\mathbb V}
\def\BBW {\mathbb W}   \def\BBX {\mathbb X}    \def\BBY {\mathbb Y}
\def\BBZ {\mathbb Z}

\def\GTA {\mathfrak A}   \def\GTB {\mathfrak B}    \def\GTC {\mathfrak C}
\def\GTD {\mathfrak D}   \def\GTE {\mathfrak E}    \def\GTF {\mathfrak F}
\def\GTG {\mathfrak G}   \def\GTH {\mathfrak H}    \def\GTI {\mathfrak I}
\def\GTJ {\mathfrak J}   \def\GTK {\mathfrak K}    \def\GTL {\mathfrak L}
\def\GTM {\mathfrak M}   \def\GTN {\mathfrak N}    \def\GTO {\mathfrak O}
\def\GTP {\mathfrak P}   \def\GTR {\mathfrak R}    \def\GTS {\mathfrak S}
\def\GTT {\mathfrak T}   \def\GTU {\mathfrak U}    \def\GTV {\mathfrak V}
\def\GTW {\mathfrak W}   \def\GTX {\mathfrak X}    \def\GTY {\mathfrak Y}
\def\GTZ {\mathfrak Z}   \def\GTQ {\mathfrak Q}

\def\sign{\mathrm{sign\,}}
\def\bdw{\prt\Gw\xspace}
\def\nabu{|\nabla u|}


\begin{abstract}
 Let $\Gw$ be a bounded domain in $\BBR^N$ and $T>0$. We study the problem 
\bel{equ} \left\{ \BA{lll} u_t - \Gd u \pm g(u) &= \gm \quad &\text{in } Q_T:=\Gw \ti (0,T) \\
\phantom{------,}
u&=0 &\text{on } \prt \Gw \times (0,T)\\
\phantom{----,}
u(.,0) &=\gw  &\text{in } \Gw. \EA \right. \ee
where $\gm$ and $\gw$ are bounded Radon measures in $Q_T$ and $\Gw$ respectively and $g(u) \sim e^{a |u|^q} $ with $a>0$ and $q \geq 1$. We provide a sufficient condition in terms of fractional maximal potentials of $\gm$ and $\gw$ for solving \eqref{equ}.
\end{abstract}


\medskip

\section{Introduction}
Let $\Gw$ be a bounded domain in $\BBR^N$ ($N \geq 2$), $T>0$ and $Q_T:=\Gw \ti (0,T)$. Denote by $\GTM^b(\Gw)$ ($\GTM^b(Q_T)$) the space of bounded Radon measures on $\Gw$ (resp. $Q_T$) and $\GTM^b_+(\Gw)$ (resp. $\GTM^b_+(Q_T)$) the positive cone of $\GTM^b(\Gw)$ (resp. $\GTM^b(Q_T)$). For $a>0$, $q \geq 1$, $\ell \geq 1$, define 
\bel{g} \CE_\ell(s)=e^s -\sum_{j=0}^{\ell-1}\frac{s^j}{j!}, s \in \BBR \text{ and } g_\ell(u)=\CE_\ell(a |u|^q). \ee  
In the present paper, we deal with the question of existence and uniqueness of solution to 
\bel{A1}\left\{ \BA{lll}
       u_t - {\Delta }u + sign(u)g_\ell(u) &= \mu \qq &\text{in } Q_T\\
       \phantom{---------,,} 
       u  &= 0 &\text{on } \prt \Gw \ti (0,T) \\ 
      \phantom{-------,,}
       u(.,0) &=\omega &\text{in }  \Gw  \\  
\EA \right.\ee    
where  $\gw \in \GTM^b(\Gw)$ and $\gm \in \GTM^b(Q_T)$. This study is inspired by recent works on elliptic equations with exponential absorption and measure data. In particular, in \cite{BLOP}, D. Bartolucci et al. proved that under the conditions $N >2$, $\gn \in \GTM^b(\Gw)$, $\gn \leq 4\pi\CH^{N-2}$ (here $\CH^{N-2}$ is $(N-2)$-dimensional Hausdorff measure in $\BBR^N$) there exists a unique solution of
\bel{A} \left\{ \BA{lll}
-\Gd u + e^u -1 &= \gn \qq &\text{in } \Gw \\
\phantom{------}
u &=0 &\text{on } \prt \Gw.
\EA \right. \ee
When $N=2$, a characterization of the set of measures for which the problem $(\ref{A})$ has a solution was given by J. L. V\'aquez (see \cite{Va}). 


Concerning the case of nonlinear operators, due to delicate estimates on {\it Wolff potentials} and {\it fractional maximal operators} (see \cite{BVH} for the definitions),  M. F. Bidaut V\'eron et al. \cite{BVHV} established a sufficient condition on $\gl \in \GTM^b(\Gw)$ for which the problem
$$ \left\{ \BA{lll}
-\Gd_p u + sign(u)g_\ell(u) &= \gl \qq &\text{in } \Gw \\
\phantom{--------,,,}
u &=0 &\text{on } \prt \Gw
\EA \right. $$
admits a {\it renormalized solution} where $\Gd_p u=div(|\nabla u|^{p-2}\nabla u)$ with $1<p<N$. 

Recently, M. F. Bidaut V\'eron and Q. H. Nguyen  have considered the parabolic problem
\bel{A2}\left\{ \BA{lll}
       u_t - \Gd_p u + sign(u)g_\ell(u) &= \mu \qq &\text{in } Q_T\\
       \phantom{----------,} 
       u  &= 0 &\text{on } \prt \Gw \ti (0,T) \\ 
      \phantom{--------,,}
       u(.,0) &=u_0 &\text{in }  \Gw \\  
\EA \right.\ee    
where $1<p<N$, $u_0 \in L^1(\Gw)$, $\gm \in \GTM^b(Q_T)$. Because of  lack of necessary tools concerning parabolic Wolff potentials, they only focused on the case where $\gm$ satisfies
\bel{ten}|\gm| \leq \gl \otimes \vartheta \ee
with $\gl \in \GTM_+^b(\Gw)$ and $\vartheta \in L^1_+((0,T))$ (here the notation $\otimes$ denotes the tensorial product). Under the condition \eqref{ten}, instead of dealing with $\gm$, they were concerned with $\gl$, which enables them to employ results developed by themselves on elliptic Wolff potentials to point out the existence of solutions to $(\ref{A2})$. 

In this paper, by limiting ourselves to the case of linear operator, we show that the condition \eqref{ten} can be removed. More precisely, when $p=2$, by adapting techniques used in \cite{BVHV} to parabolic framework, we obtain a sufficient condition on $\gm \in \GTM^b(Q_T)$ and $u_0 \in \GTM^b(\Gw)$ respectively in terms of parabolic and elliptic fractional maximal operators for solvability of \eqref{A2} . In order to state the results, we first introduce some notations. \smallskip

\noindent{\bf Notations and terminology.} For $\ga > 0$ and $\gb \geq 0$, set 
$$ h_{1,\ga}(s)= (-\ln (s\wedge 2^{-1}))^{\frac{1}{\ga}}, \q h_{2,\gb}(s)= \left(\ln \left(2ds^{-1} \lor 2 \right) \right)^{-\gb},\forevery s>0$$
where $d=diam(\Gw)+T$ (here $a \wedge b = \min\{a,b\}$, $a \lor b =\max\{a,b \}$). For $0<R\leq \infty$, we denote the {\it $R-$truncated $\ga-$fractional maximal potential of $\gw$} by
$$ \BBM_{\ga,R}^{1}[\gw](x)=\sup_{0 < s \leq R}\left( \myfrac{\gw(B_s(x))}{s^Nh_{1,\ga}(s)} \right) \q \text{for a.e. } x \in \BBR^N$$ 
where $B_s(x)$ is the ball of center $x$ and radius $s>0$. The {\it parabolic $R-$truncated $\gb-$fractional maximal potential of $\gm$} is defined by
$$ \BBM_{\gb,R}^{2}[\gm](x,t)=\sup_{0 < s \leq R}\left( \myfrac{\gm(Q_s(x,t))}{s^Nh_{2,\gb}(s)} \right) \q \text{for a.e.} (x,t) \in \BBR^{N+1}$$
where $Q_s(x,t)=B_s(x)\ti(t-s^2/2,t+s^2/2)$. Finally, the {\it parabolic $R-$truncated Wolff potential of $\gm$} is defined by
$$ \BBW_R[\gm](x,t)=\myint{0}{R}\myfrac{\gm(Q_s(x,t))}{s^N}\myfrac{ds}{s} \q \text{for a.e.} (x,t) \in \BBR^{N+1}.$$
\bdef{sol} Let $f \in C(\BBR)$, $\gm \in \GTM^b(Q_T)$ and $\gw \in \GTM^b(\Gw)$. A function $u$ is a solution of 
\bel{geneq}\left\{ \BA{lll}
       u_t -  \Gd u + f(u)&= \mu \qq &\text{in } Q_T \\ 
      \phantom{-----,,,}
       u  &= 0 &\text{on } \prt \Gw \ti (0,T) \\ 
      \phantom{----,}
       u(.,0) &=\omega &\text{in }  \Gw  \\  
\EA \right.\ee   
if $u \in L^1(Q_T)$, $f(u) \in L^1(Q_T)$ and 
\bel{sol} \myint{Q_T}{}(-u(\gz_t + \Gd \gz) + f(u)\gz)dxdt = \myint{Q_T}{}\gz d\gm + \myint{\Gw}{}\gz(.,0)d\gw\ee
for every $\gz \in X(Q_T)$, which is the space of functions in $C^{2,1}(\ovl Q_T)$ vanishing on $(\prt \Gw \ti [0,T]) \cup (\Gw \ti \{T\})$.  
\es
In the sequel, if $\gm \in \GTM(Q_T)$ ($\gw \in \GTM(\Gw)$ resp.), we will consider $\gm$  (resp. $\gw$) as a measure in $\BBR^{N+1}$  (resp. $\BBR^N$) vanishing outside of $Q_T$ (resp. $\Gw$). 
The first result in the paper is the following
\bth{exist2} Let $\Gw$ be a bounded domain with $C^2$ boundary. Assume $a>0$, $q \geq 1$, $\ell \geq 1$, $\ga \geq q$, $\gb \in [\frac{q-1}{q},1)$, $f_1 \in L^1(\Gw)$, $f_2 \in L^1(Q_T)$, $\gw \in \GTM^b(\Gw)$ and $\gm \in \GTM^b(Q_T)$. There exist $M_1=M_1(N,\ga,a)$ and $M_2=M_2(N,\gb,a)$ such that if $\norm{\BBM_{\ga,\infty}^{1}[\gw^{\pm}]}_{L^\infty(\BBR^N)} < M_1$ and $\norm{\BBM_{\gb,\infty}^{2}[\gm^{\pm}]}_{L^\infty(\BBR^{N+1})}<M_2$ then  the problem
\bel{absorp-f}\left\{ \BA{lll}
       u_t -  \Gd u + sign(u)g_\ell(u)&= \mu + f_2\qq &\text{in } Q_T \\ 
      \phantom{---------,,}
       u  &= 0 &\text{on } \prt \Gw \ti (0,T) \\ 
      \phantom{-------,,}
       u(.,0) &=\omega + f_1&\text{in }  \Gw  \\  
\EA \right.\ee   
admits a unique solution $u$ satisfying $e^{a|u|^q} \in L^1(Q_T)$.
\es
We also consider the problem associated to equation with source terms
\bel{source}\left\{ \BA{lll}
       u_t -  \Gd u &= g_\ell(u) + \mu \qq &\text{in } Q_T \\ 
      \phantom{--,,,}
       u  &= 0 &\text{on } \prt \Gw \ti (0,T) \\ 
      \phantom{-,}
       u(.,0) &=\omega &\text{in }  \Gw  \\  
\EA \right.\ee   
where $\gw \in \GTM_+^b(\Gw)$, $\gm \in \GTM_+^b(Q_T)$ and $g_\ell$ is defined as in $(\ref{g})$ $a>0$, $p \geq 1$, $ \ell p>1$. \medskip

Let $G(x,t)$ be the heat kernel in $\BBR^N$ which is defined by $G(x,t)=(4\pi t)^{-\frac{N}{2}}e^{-\frac{\abs{x}^2}{4t}}$  if $x \in \BBR^N, t>0$ and $G(x,t)=0$ if $x \in \BBR^N, t \geq 0$. For any $y \in \Gw$, denote by $G^\Gw(x,t,y)$ the fundamental solution of the heat equation in $\Gw$ with zero Dirichlet condition on $\prt \Gw$ and initial condition $\gd_y$ (Dirac measure concentrated at $y$). Clearly
$ G^\Gw(x,t,y) \leq G(x-y,t)$ for every $x,y \in \Gw, t>0$. If $\gw \in \GTM^b(\Gw)$, we denote $\BBG[\gw](x,t)=\int_{\Gw}G^\Gw(x,y,t)d\gw(y)$. \medskip

Existence result for \eqref{source} is stated in the following theorem
\bth{exist3} Let $\Gw$ be a bounded domain with $C^2$ boundary. Assume $a>0$, $q \geq 1$, $\ell q>1$, $\ga \geq q$, $\gb \geq [\frac{q-1}{q},1)$, $\gw \in \GTM_+^b(\Gw)$ and  $\gm \in \GTM_+^b(Q_T)$.  There exist $c=c(N)$,
$b_0=b_0(N,d,\ell,q) \in (0,1]$, $M_1=M_1(N,a,\ga,q,\ell,d)$ and $M_2=M_2(N,a,\gb,q,\ell,d)$ such that if $\norm{\BBM_{\ga,\infty}^{1}[\gw]}_{L^\infty(\BBR^N)} \leq M_1$ and $\norm{\BBM_{\gb,\infty}^2[\gm]}_{L^\infty(\BBR^{N+1})} \leq M_2$ then the problem $(\ref{source})$ admits a nonnegative solution $u$ which satisfies
\bel{estu} u \leq \BBG[\gw] + c \BBW_{2d}[\gm] + c b_0 . \ee
\bel{exploc} g_\ell(2\BBG[\gw]+2c \BBW_{2d}[\gm]+2c b_0) \in L^1(Q_T). \ee
\es \medskip

The paper is organized as follows. In Section 2, we establish estimates on parabolic Wolff potentials. Section 3 is devoted to the study of linear parabolic equations with measure data. In section 4 we apply results obtained in Section 2 and Section 3 to prove existence of solution to equation $(\ref{absorp})$ and  $(\ref{source})$.  \medskip

\noindent{\bf Acknowledgements.} The author warmly thanks Q. H. Nguyen for fruitful discussions. He is grateful to the editor and the referee for their useful comments. The author is supported by the Israel Science Foundation, grant No. 91/10, founded by the Israel Academy of Sciences and Humanities.
\section{Estimates on parabolic Wolff potentials}
We start this section with some notations. If $A$ is a measurable set in $\BBR^{N+1}$, we denote by $|A|$ the Lebesgue measure of $A$. If $f$, $g$ are functions defined in $\BBR^{N+1}$ and $a,b \in \BBR$ then set $\{f >a\}:=\{(x,t) \in \BBR^{N+1}: f(x,t) >a \}$, $\{ f>a, g \leq b\}:=\{f>a\} \cap \{g\leq b\}$. Finally $\chi_A$ denote the characteristic function of $A$.
\bprop{estset} Assume $\gb \in [0,1)$ and $r>0$. There exist $c_1=c_1(N,\gb)$ and $\ge_1=\ge_1(N,\gb,d,r)$ such that, for any $\gm \in \GTM_+(\BBR^{N+1})$ satisfying $supp(\gm) \sbs Q_r(x^*,t^*)$ for some $(x^*,t^*) \in \BBR^N \ti \BBR$ and for any $R \in (0,\infty]$, $\ge \in (0,\ge_1]$, $\gl > \gm(\BBR^{N+1})l(r,R)$ there holds
\bel{estset} \BA{lll} \abs{\{\BBW_R[\gm]>3\gl, \BBM_{\gb,R}^{2}[\gm] \leq \ge \gl \}} \\
\phantom{---}
\leq c_1\exp\left(-2^{-\frac{\gb}{1-\gb}}(1-\gb)^{\frac{1}{1-\gb}}\ge^{-\frac{1}{1-\gb}}\ln2\right)\abs{\{\BBW_R[\gm]>\gl\}}
\EA \ee
where $l(r,R)=N^{-1}((r \wedge R)^{-N} -R^{-N})$ if $R<\infty$ and $l(r,\infty)=N^{-1}r^{-N}$. If $\gb=0$ then $\ge_0$ depends only on $N$, $\gb$ and $(\ref{estset})$ holds true for every $\gm \in \GTM_+(\BBR^{N+1})$ with compact support in $\BBR^{N+1}$, $R \in (0,\infty]$, $\ge \in (0,\ge_1]$, $\gl>0$.
\es
\Proof We adapt the ideas used in \cite{BVHV} to parabolic setting. Denote the parabolic distance by 
$$ d_P((x,t),(y,\gt))=\abs{x-y}+\abs{t-\gt}^{1/2} \quad \forall (x,t), (y,\gt) \in \BBR^N \ti \BBR.$$ 
If $A,B \sbs \BBR^{N+1}$, we denote 
$$diam(A)=\sup\{d_P((x,t),(y,\gt)): (x,t), (y,\gt) \in A\}, $$ 
$$\dist(A,B)=inf\{d_P((x,t),(y,\gt)): (x,t) \in A, (y,\gt) \in B\}.$$ 
For any $x=(x^1,..,x^N) \in \BBR^N$, $t \in \BBR$, $r>0$, the parabolic cube of center $x$ and edge $r$ is defined as follows
$$ K_r(x,t)=\left[x^1-\frac{r}{2},x^1+\frac{r}{2}\right]\ti ...\ti \left[x^N-\frac{r}{2},x^N+\frac{r}{2}\right]\ti \left[t-\frac{r^2}{2},t+\frac{r^2}{2}\right].$$
Notice that $diam(K_r(x,t))=(\sqrt{N}+1)r$ for every $(x,t) \in \BBR^N \ti \BBR$. \medskip

\noindent{\bf Case 1: $R=\infty$.} Let $\gl>0$ and  $D_\gl=\{\BBW_\infty[\gm]>\gl\}$. By Whitney covering lemma (see \cite{LM}), there exists a countable family $\CK:=\{K_i\}$, where $K_i=K_{r_i}(x_i,t_i)$, such that $\cup_i K_i =D_\gl$, $\mathring{K_i} \cap \mathring{K_j}= \emptyset$ if $i \neq j$ and there exists a positive constant $C_w=C_w(N)>1$ such that
$$C_w^{-1}diam(K_i) \leq \dist(K_i,D_\gl^c) \leq C_w diam(K_i) \forevery i. $$
Let $\ge>0$ and  denote $F_{\ge,\gl}=\{\BBW_\infty[\gm]>3\gl, \BBM_{\gb,\infty}^{2}[\gm] \leq \ge\gl\}$. 
We will show that there exist $c_2=c_2(N,\gb)>0$ and $\ge_1=\ge_1(N,\gb,r,d)$ such that for any $K \in \CK$, $\ge \in (0,\ge_1]$ and $\gl>(\gm(\BBR^{N+1}))l(r,\infty)$ there holds
\bel{q1} |F_{\ge,\gl}\cap K| \leq c_2 \exp\left(-2^{-\frac{\gb}{1-\gb}}(1-\gb)^{\frac{1}{1-\gb}}\ge^{-\frac{1}{1-\gb}}\ln2\right) |K|. \ee
In order to do that we prove \medskip

\noindent{\bf Assertion 1:} There exists $\ge_2=\ge_2(N,\gb)$ such that for any $K \in \CK$, $\ge \in (0,\ge_2]$, $\gl>0$, there holds $ F_{\ge,\gl} \cap K \sbs E_{\ge,\gl} $ where 
$$E_{\ge,\gl}=\{ (x,t) \in K: \BBW_{(1+C_w)diam(K)}[\gm](x,t) > \gl, \BBM_{\gb,\infty}^{2}[\gm](x,t) \leq \ge \gl \}.$$

Take $K \in \CK$ such that $ F_{\ge,\gl} \cap K \ne \emptyset$ and take $({\tl x},{\tl t}) \in D_\gl^c$ satisfying $\dist(({\tl x},{\tl t}),K) \leq C_w diam(K)$. Denote $r_0=(1+C_w)diam(K)$. For any $k \in \BBN$ and $(x,t) \in F_{\ge,\gl} \cap K$, we denote
$$ A_k=\myint{2^kr_0}{2^k\frac{1+2^{k+1}}{1+2^k}r_0}\myfrac{\gm(Q_s(x,t))}{s^N}\myfrac{ds}{s}, \q B_k=\myint{2^k\frac{1+2^{k+1}}{1+2^k}r_0}{2^{k+1}r_0}\myfrac{\gm(Q_s(x,t))}{s^N}\myfrac{ds}{s}.$$
Note that $ B_k \leq c_3\ge\gl2^{-k}$ where $c_3=c_3(\gb)$. Set $\gd=\left( \frac{2^k}{1+2^k} \right)^N$ then $1-\gd < c_4 2^{-k}$ with $c_4=c_4(N)$. Consequently, $(1-\gd)A_k \leq c_5\ge \gl 2^{-k}$ with $c_5=c_5(N,\gb)$. For any $(x,t) \in F_{\ge,\gl} \cap K$ and $s \in [(1+2^k)r_0, (1+2^{k+1})r_0]$, we have $ Q_{\frac{2^k}{1+2^k}s}(x,t) \sbs Q_s({\tl x},{\tl t})$, from which it follows
$$ \gd A_k \leq \myint{(1+2^k)r_0}{(1+2^{k+1})r_0}\myfrac{\gm(Q_s({\tl x},{\tl t}))}{s^N}\myfrac{ds}{s}. $$
As a consequence,
$$ \myint{2^kr_0}{2^{k+1}r_0}\myfrac{\gm(Q_s(x,t))}{s^N}\myfrac{ds}{s}=A_k+B_k\leq c_6 2^{-k}\ge \gl +  \myint{(1+2^k)r_0}{(1+2^{k+1})r_0}\myfrac{\gm(Q_s({\tl x},{\tl t}))}{s^N}\myfrac{ds}{s}$$
where $c_6=c_6(N,\gb)$. Therefore
$$  \myint{r_0}{\infty}\myfrac{\gm(Q_s(x,t))}{s^N}\myfrac{ds}{s} \leq 2c_6\ge \gl +  \myint{2r_0}{\infty}\myfrac{\gm(Q_s({\tl x},{\tl t}))}{s^N}\myfrac{ds}{s} \leq (1+2c_6\ge)\gl.$$
Put $\ge_2=(2c_6)^{-1}$. If $\ge \in (0,\ge_2]$ then 
$$ \myint{r_0}{\infty}\myfrac{\gm(Q_s(x,t))}{s^N}\myfrac{ds}{s} \leq 2\gl, $$
which implies Assertion 1. \medskip

\noindent{\bf Assertion 2:} There exists $c_7=c_7(N,\gb)$ such that
\bel{estE} \abs{E_{\ge,\gl}} \leq c_7\exp\left(-2^{-\frac{\gb}{1-\gb}}\ln 2 (1-\gb)^{\frac{1}{1-\gb}}\ge^{-\frac{1}{1-\gb}}\right)\abs{K}. \ee
Denote $Q^*_1:=Q_r(x^*,t^*)$ and $Q^*_2:=Q_{2r}(x^*,t^*)$ for some $r>0$ and $(x^*,t^*) \in \BBR^{N+1}$. Let $\gl>(\gm(\BBR^{N+1})l(r,\infty)$. If $(x,t) \in (Q^*_2)^c$ and $s<r$ then $Q_s(x,t) \cap Q^*_1=\emptyset$.
Hence
$$ \BBW_\infty[\gm](x,t)=\myint{r}{\infty}\myfrac{\gm(Q_s(x,t))}{s^N}\myfrac{ds}{s} \leq \gm(\BBR^{N+1})l(r,\infty). $$
Therefore $D_\gl \sbs Q^*_2$, which in turn implies $r_0 \leq 5(1+C_w)r$. Next, we set $m_0=(\ln 2)^{-1}\max(1,\ln (5d^{-1}(1+C_w)r))$ then $2^{-m}r_0 \leq d$ for every $m \geq m_0$. For any $(x,t) \in E_{\ge,\gl}$ and $ m>m_0^{\frac{1}{1-\gb}}$,
$$ \BA{lll} \myint{2^{-m}r_0}{r_0}\myfrac{\gm(Q_s(x,t))}{s^N}\myfrac{ds}{s} \leq \myfrac{\ge\gl}{1-\gb}((m-m_0)\ln  2)^{1-\gb} + m_0\ge\gl \leq \myfrac{2\ge\gl}{1-\gb}m^{1-\gb}. \EA $$ 
If we define 
$$h_i(x,t)=\myint{2^{-i}r_0}{2^{-i+1}r_0}\myfrac{\gm(Q_s(x,t))}{s^N}\myfrac{ds}{s}, \quad i \in \BBN,$$
then for any $m \geq m_0^{\frac{1}{1-\gb}}$,
$$ \BBW_{r_0}[\gm](x,t) \leq \myfrac{2\ge\gl}{1-\gb}m^{1-\gb} + \mysum{i=m+1}{\infty}h_i(x,t). $$
Consequently, for $0<\gg<2$,
$$ \BA{lll} |E_{\ge,\gl}| \\ 
\phantom{} \leq \mysum{i=m+1}{\infty}\abs{\left\{ (x,t) \in K: h_i(x,t) > 2^{-\gg(i-m-1)}(1-2^{-\gg})(1-\frac{2}{1-\gb}m^{1-\gb}\ge)\gl    \right\}}. \EA $$
After a long computation, we get
$$ \BA{lll} \abs{\{ (x,t) \in K: h_i(x,t)>s \}}  \leq \myfrac{2^{2+2N}(\ln 2)^{-\gb}}{s}2^{-2i}r_0^{N+2}\ge\gl  
\leq c_8\myfrac{2^{-2i}}{s}|K|\ge\gl \forevery i
\EA $$
which leads to
\bel{ul1} \abs{E_{\ge,\gl}} \leq c_9 2^{-2m}\myfrac{\ge}{1-\frac{2}{1-\gb}m^{1-\gb}\ge}|K| \forevery m>m_0^{\frac{1}{1-\gb}} \ee
where $c_9=c_9(N,\gb)$. Set $\ge_1=\min\{\frac{1}{4(1-\gb)^{-1}m_0} ,\ge_2\}$. For any $\ge \in (0,\ge_1]$, we choose $m \in \BBN$ such that
$$ \left( \myfrac{1-\gb}{2} \right)^{\frac{1}{1-\gb}} \left( \myfrac{1}{\ge}-1 \right)^{\frac{1}{1-\gb}}-1 <m \leq  \left( \myfrac{1-\gb}{2} \right)^{\frac{1}{1-\gb}} \left( \myfrac{1}{\ge}-1 \right)^{\frac{1}{1-\gb}}. $$
Then \bel{ul2} \myfrac{\ge}{1-\frac{2}{1-\gb}m^{1-\gb}\ge} \leq 1 \text{~~~and~~~}
2^{-2m} \leq 4\exp\left(-2^{-\frac{\gb}{1-\gb}}\ln 2 (1-\gb)^{\frac{1}{1-\gb}}\ge^{-\frac{1}{1-\gb}}\right).\ee
Combining $(\ref{ul1})$-$(\ref{ul2})$ yields to Assertion 2. Finally $(\ref{q1})$ follows straightforward.

If $\gb=0$ then for any $m \in \BBN$, $\ge>0$, $\gl>0$ and $(x,t) \in E_{\ge,\gl}$
$$ \BBW_{r_0}[\gm](x,t) \leq \ge\gl m + \mysum{i=m+1}{\infty}h_i(x,t).$$
Consequently, with $m\ge<1$, $ \abs{E_{\ge,\gl}} \leq c_92^{-2m}\ge(1-m\ge)^{-1}|K|$.
Put $\ge_1=\min\{\frac{1}{2},\ge_2\}$ then for any $\ge \in (0,\ge_1]$ and $\ge^{-1}-2 < m <\ge^{-1}-1$, we obtain
$$ \abs{E_{\ge,\gl}} \leq 16c_9\exp(-2\ge^{-1}\ln2 )|K|, $$
which leads to $(\ref{q1})$. \medskip

\noindent{\bf Case 2: $R<\infty$.} For $\gl>0$, $D^R_{\gl}=\{\BBW_R[\gm]>\gl\}$ is an open subset of $\BBR^{N+1}$. By Whitney covering lemma, there exists a countable family of closed cubes $\CK:=\{K_i\}$ such that $\cup_i K_i =D^R_{\gl}$, $\mathring{K_i} \cap \mathring{K_j}= \emptyset$ if $i \neq j$ and $\dist(K_i,(D^R_{\gl})^c) \leq C_w diam(K_i)$. If $K \in \CK$ is such that $diam(K)>\frac{R}{2C_w}$, there exists a finite number $n_K$ of closed dyadic cubes $\{P_{j,K}\}_{j=1}^{n_K}$ satisfying $\cup_{j=1}^{n_K}P_{j,K}=K$, $\mathring{P_{i,K}}\cup \mathring{P_{j,K}}=\emptyset$ if $i \neq j$ and $\frac{R}{4C_w}<diam(P_{j,K})<\frac{R}{2C_w}$. We set $\CK'=\{K \in \CK: diam(K) \leq \frac{R}{2C_w}\}$, $\CK''=\{P_{i,K}:1 \leq i \leq n_K, K \in \CK, diam(K) > \frac{R}{2C_w}\}$ and $\tl \CK=\CK' \cup \CK''$.

For $\ge>0$, we denote $F^R_{\ge,\gl}=\{\BBW_R[\gm]>3\gl, \BBM_{\gb,R}^{2}[\gm] \leq \ge\gl\}$.
Let $K \in \tl \CK$ such that $F^R_{\ge,\gl} \cap K \neq \emptyset$ and set $r_0=(1+C_w)diam(K)$. \medskip

\noindent{\bf Case 2.i: $\dist((D^R_\gl)^c,K) \leq C_w diam(K) $}. Let $({\tl x},{\tl t}) \in (D^R_\gl)^c$ such that $\dist(({\tl x},{\tl t}),K) \leq C_wdiam(K)$ and $\BBW_R[\gm]({\tl x},{\tl t}) \leq \gl$.  By using the same argument as in Case 1, we deduce that for any $(x,t) \in K \cap F^R_{\ge,\gl}$,
$$ \myint{r_0}{R}\myfrac{\gm(Q_s(x,t))}{s^N}\myfrac{ds}{s} \leq (1 + c_{10}\ge)\gl $$
where $c_{10}=c_{10}(N,\gb)$. \medskip

\noindent{\bf Case 2.ii: $\dist((D^R_\gl)^c,K) > C_w diam(K) $}. Then $K \in \CK''$ and hence $\frac{R}{4C_w} < diam(K) \leq \frac{R}{2C_w}$. Therefore, for any $(x,t) \in K \cap F^R_{\ge,\gl}$,
$$ \BA{lll} \myint{r_0}{R}\myfrac{\gm(Q_s(x,t))}{s^N}\myfrac{ds}{s} &\leq \myint{\frac{1+C_w}{4C_w}R}{R}\myfrac{\gm(Q_s(x,t))}{s^N}\myfrac{ds}{s} \\
&= \ge\gl (\ln 2)^{-\gb} \ln\left(\frac{4C_w}{1+C_w}\right) \leq 2\ge\gl. \EA $$

Put $\ge_3=\min\{1,c_{10}^{-1}\}$ then for any $\ge \in (0,\ge_3]$, $K \cap F^R_{\ge,\gl} \sbs E^R_{\ge,\gl}$ where 
$$E^R_{\ge,\gl}=\{ (x,t) \in K: \BBW_{r_0}[\gm](x,t) > \gl, \BBM_{\gb,R}^{2}[\gm](x,t) \leq \ge \gl \}.$$
By proceeding as in case 1, we can derive \eqref{estE} with $E_{\ge,\gl}$ replaced by $E^R_{\ge,\gl}$ and with another constant. Thus \eqref{q1} follows straightforward. Finally, the case $\gb=0$ is treated as in case 1.  \qed

\bth{double} Assume $0 \leq \gb <1$ and $r>0$. Set $\gd_1=2(\frac{1-\gb}{6})^{\frac{1}{1-\gb}}\ln2$. There exists $c_{11}=c_{11}(N,\gb,d,r)$ such that for any $R \in (0,\infty]$, $\gd \in (0,\gd_1)$, $\gm \in \GTM_+(\BBR^{N+1})$, $r' \in (0,r]$, $(x^*,t^*) \in \BBR^N \ti \BBR$, there holds
\bel{double} \myfrac{1}{\abs{Q_{2r'}(x^*,t^*)}}\myint{Q_{2r'}(x^*,t^*)}{}\exp\left( \gd M_2^{-\frac{1}{1-\gb}}\BBW_R[\gm^*]^{\frac{1}{1-\gb}}\right)dxdt \leq \myfrac{c_{11}}{\gd_1-\gd} \ee
where $\gm^*=\gm \chi_{Q_{r'}(x^*,t^*)}$ and $M_2=\norm{\BBM_{\gb,R}^{2}[\gm^*]}_{L^\infty(\BBR^{N+1})}$. If $\gb=0$ then $c_{11}$ is independent of $r$.
\es
\Proof Let $\gm \in \GTM_+(\BBR^{N+1})$ satisfy $M_2<\infty$. Due to \rprop{estset} there exist $c_{1}=c_{1}(N,\gb)$ and $\ge_1=\ge_1(N,d,\gb,r)$ such that for any $R \in (0,\infty]$, $\ge \in (0,\ge_1]$, $\gl>\gm(\BBR^{N+1})l(r',R)$ there holds
\bel{estset11} \BA{lll} \abs{\{\BBW_R[\gm^*]>3\gl, \BBM_{\gb,R}^{2}[\gm^*] \leq \ge \gl \}} \\
\phantom{----}
\leq c_{1}\exp\left(-2^{-\frac{\gb}{1-\gb}}(1-\gb)^{\frac{1}{1-\gb}}\ge^{-\frac{1}{1-\gb}}\ln2\right)\abs{\{\BBW_R[\gm^*]>\gl\}}. \EA \ee
Since $ \gm^*(\BBR^{N+1})l(r',R) <N^{-1}(\ln 2)^{-\gb}M_2 $, we can choose $\ge$ and $\gl$ in $(\ref{estset11})$ such that $ \ge=\gl^{-1}M_2$ with $\gl>\max\{\ge_1^{-1}, N^{-1}(\ln 2)^{-\gb}\}M_2$.
By using similar argument as in \rprop{estset}, we deduce that $D_\gl^* \sbs Q'_2$ where $D_\gl^*=\{\BBW_R[\gm^*]>\gl\}$ and $Q'_2=Q_{2r'}(x^*,t^*)$. Hence
\bel{estset12} \BA{lll} \abs{\{\BBW_R[\gm^*]>3\gl\} \cap Q'_2} \\ 
\phantom{---}
\leq c_{1}\exp\left(-2^{-\frac{\gb}{1-\gb}}(1-\gb)^{\frac{1}{1-\gb}}M_2^{-\frac{1}{1-\gb}}\ln2\gl^{\frac{1}{1-\gb}}\right)\abs{Q'_2}\EA \ee
Therefore
$ \abs{\{ \Psi > \gth\} \cap Q'_2 } \leq |Q'_2|\chi_{(0,\gth_0]} + c_{1}e^{-\gd_1\gth}|Q'_2|\chi_{(\gth_0,\infty)} $
where $\Psi=M_2^{-\frac{1}{1-\gb}}\BBW_R[\gm^*]^{\frac{1}{1-\gb}}$ and $\gth_0=(3\max\{\ge_1^{-1},N^{-1}(\ln2)^{-\gb}\})^{\frac{1}{1-\gb}}$. Thus, for each $\gd \in (0,\gd_1)$,
$$ \myint{Q'_2}{}e^{\gd \Psi}dxdt \leq (e^{\gd \gth_0}-1)|Q'_2| + \myfrac{c_{1}\gd}{\gd_1-\gd}|Q'_2| $$
which implies the desired estimate. \qed

The next result is crucial for proving existence of solution to \eqref{absorp} and \eqref{source} in section 4. 
\bth{Hexp} Assume that $0 \leq \gb <1$, $R>0$ and $\gm \in \GTM_+(\BBR^{N+1})$ satisfies  $\norm{\BBM_{\gb,\infty}^{2}[\gm]}_{L^\infty(\BBR^{N+1})} \leq M_2$. Then there exist $\gd_2=\gd_2(\gb)$ and $c_i=c_i(N,\gb,R,d)$ $(i=12,13)$ such that for any $r \in (0,R)$ and any $(x,t) \in \BBR^{N+1}$, there holds
\bel{expH} \myint{Q_r(x,t)}{}\exp(\gd_2 M_2^{-\frac{1}{1-\gb}} (\BBW_R[\gm])^{\frac{1}{1-\gb}})dy d\tau<c_{12},\ee
\bel{Hexp} \norm{\BBW_R[\exp(\gd_2M_2^{-\frac{1}{1-\gb}}(\BBW_R[\gm])^{\frac{1}{1-\gb}})]}_{L^\infty(\BBR^{N+1})} < c_{13}. \ee
\es
\Proof Fix $(x,t) \in \BBR^N \ti \BBR$. For every $(y,\gt) \in \BBR^N \ti \BBR$, we have
$$ \BA{lll}\BBW_R[\gm](y,\gt)=\BBW_r[\gm](y,\gt) + \myint{r}{R}\myfrac{\gm(Q_s(y,\gt))}{s^N}\myfrac{ds}{s}\\
\phantom{}
\leq \BBW_r[\gm](y,\gt)+M_2\myint{r \wedge d}{d}(\ln(2ds^{-1}))^{-\gb}\myfrac{ds}{s} + M_2\myint{d}{R \lor d}(\ln2)^{-\gb}\myfrac{ds}{s}\\
\phantom{}
\leq  \BBW_r[\gm](y,\gt)+M_2(\ln2)^{-\gb}\ln\left(\myfrac{R}{d} \lor 1\right) + M_2(1-\gb)^{-1}\left(\ln\left(\myfrac{d}{r} \lor 1\right)\right)^{1-\gb}.
\EA $$
Consequently,
$$ \BA{lll}\BBW_R[\gm](y,\gt)^{\frac{1}{1-\gb}} \leq 3^{\frac{\gb}{1-\gb}}M_2^{\frac{1}{1-\gb}}(\ln2)^{-\frac{\gb}{1-\gb}}\left( \ln\left( \myfrac{R}{d} \lor 1 \right)\right)^{\frac{1}{1-\gb}} \\
\phantom{-----}
+3^{\frac{\gb}{1-\gb}}\BBW_r[\gm](y,\gt)^{\frac{1}{1-\gb}} + 3^{\frac{\gb}{1-\gb}}M_2^{\frac{1}{1-\gb}}(1-\gb)^{-\frac{1}{1-\gb}}\ln\left( \myfrac{d}{r} \lor 1\right).
\EA $$
Let $\gk \in (0,1]$ (to be made precise later on). It follows from the above estimate that
$$ \BA{lll} \exp\left( \myfrac{\gk \gd_1}{4.3^{\frac{\gb}{1-\gb}}}M_2^{-\frac{1}{1-\gb}}\BBW_R[\gm]^{\frac{1}{1-\gb}} \right) \\
\phantom{-----} 
\leq 2^{-1}\exp\left(\myfrac{ \gd_1}{2} M_2^{-\frac{1}{1-\gb}}\BBW_r[\gm]^{\frac{1}{1-\gb}}\right)+c_{14}\left( \myfrac{d}{r} \lor 1 \right)^{\gk c_{15}} \EA $$
where $c_{14}=c_{14}(\gb,R,d)$, $c_{15}=c_{15}(\gb)$ and $\gd_1$ is defined in \rth{double}.   For every $s \in (0,r]$, $(y,\gt) \in Q_r(x,t)$, we get $Q_s(y,\gt) \sbs Q_{2r}(x,t)$. Therefore $\BBW_r[\gm](y,\gt)=\BBW_r[\gm\chi_{Q_{2r}(x,t)}](y,\gt)$ for every $(y,\gt) \in Q_r(x,t)$. Thanks to \rth{double}, we get
$$ \myint{Q_r(x,t)}{}\exp\left(\myfrac{\gd_1}{2} M_2^{-\frac{1}{1-\gb}}(\BBW_r[\gm])^{\frac{1}{1-\gb}}\right)dyd\gt \leq c_{16}r^{N+2} $$
where $c_{16}=c_{16}(N,\gb,d,R)$. Thus
\bel{ab} \BA{lll} \myint{Q_r(x,t)}{}\exp\left( \myfrac{\gk  \gd_1}{4.3^{\frac{\gb}{1-\gb}}}M_2^{-\frac{1}{1-\gb}}\BBW_R[\gm]^{\frac{1}{1-\gb}} \right)dyd\gt \\ 
\phantom{---------}
\leq 2^{-1}c_{16}r^{N+2}+c_{17}\left(\myfrac{d}{r} \lor 1 \right)^{\gk c_{15}}r^{N+2} \EA \ee
with $c_{17}=c_{17}(N,\gb,R,d)$. By taking $\gk=1 \wedge (2c_{15})^{-1}$ and $\gd_2=2^{-2}3^{-\frac{\gb}{1-\gb}}\gk \gd_1 $, we derive $(\ref{expH})$. Finally, $(\ref{Hexp})$ follows from $(\ref{ab})$. \qed
%

\section{Estimates on solutions to linear parabolic equation}
In this section, let $\Gw$ be a bounded domain with $C^2$ boundary. We first give some estimates on solutions to 
homogeneous linear equations with initial measure data.
\bth{nononlinear0}
Assume $\ga \geq 1$, $\gd>0$ and $\omega\in \GTM^b_+(\Omega)$. There exists a positive constant $M_1=M_1(N,\ga,\gd)$ such that if  $\norm{\BBM_{\ga,\infty}^{1}[\gw]}_{L^\infty(\BBR^N)} \leq M_1$  then the unique solution $u$ to the problem 
\bel{nononlinear0}\left\{ \BA{lll}
       u_t - {\Delta }u &= 0 \qq &\text{in } Q_T \\ 
      \phantom{---}
       u  &= 0 &\text{on } \prt \Gw \ti (0,T) \\ 
      \phantom{-}
       u(.,0) &=\omega &\text{in }  \Gw . \\  
\EA \right.\ee  
satisfies 
\bel{eu} \exp ( {\delta {u^\ga }(x,t)} ) \leq c_{18}t^{-\frac{1}{2}} +2 \quad \forall (x,t) \in Q_T \ee
where $c_{18}=c_{18}(N)$.
\es 
\Proof The unique solution of $(\ref{nononlinear0})$ is represented by (see \cite{V2})
\bel{repsol0} u(x,t) = \BBG[\gw](x,t) \leq \myint{\mathbb{R}^N}{}{(4\pi t)^{-\frac{N}{2}}}e^{-\frac{|x-y|^2}{4t}}d\omega(y)  \quad \forall (x,t) \in Q_T.\ee
Fix $(x,t) \in Q_T$. Using Fubini Theorem we get 
 $$ \BA{lll}
   \myint{\mathbb{R}^N}{} (4\pi t)^{-\frac{N}{2}}e^{ { - \frac{{|x - y{|^2}}}{{4t}}} }d\omega (y) & = \myint{\mathbb{R}^N}{} (4\pi t)^{-\frac{N}{2}}\myint{\frac{{|x - y{|^2}}}{{4t}}}{\infty}  {e^{-r}dr} d\omega (y)  \\[3mm] 
    &= \myint{\mathbb{R}^N}{} {\myint{0}{\infty}  (4\pi t)^{-\frac{N}{2}}{\chi _{B_{\sqrt {4tr}}(x)}}(y)e^{-r}drd\omega (y)}   \\[3mm]
&= \myint{0}{\infty}  (4\pi t)^{-\frac{N}{2}}\omega ( B_{\sqrt {4tr}}(x))e^{- r}dr 
  \EA $$
Let $M_1>0$ be made precise later on. If  $\norm{\BBM_{\ga,\infty}^{1}[\gw]}_{L^\infty(\BBR^N)} \leq M_1$ then by combining the assumption and $(\ref{repsol0})$, we get
  \bel{e-ini1} \BA{lll}
  u(x,t) &\leq M_1\myint{0}{\infty}  (4\pi t)^{-\frac{N}{2}}(4tr)^{\frac{N}{2}}{{\left( { - \ln \left( { (4tr)^{\frac{1}{2}}  \wedge 2^{-1}} \right)} \right)}^{\frac{1}{\ga }}}e^{-r}dr  \\ [3mm]
      &= M_1\myint{0}{\infty} {{\pi ^{ - \frac{N}{2}}}{{\left( { - \ln \left( { (4tr)^{\frac{1}{2}}  \wedge 2^{-1}} \right)} \right)}^{\frac{1}{\ga }}}{r^{\frac{N}{2}}}e^{-r}dr}  \\ [3mm]
      &= M_1\myint{0}{\infty}  {{c_{19}}{{\left( { - \ln \left( { (4tr)^{\frac{1}{2}}  \wedge 2^{-1}} \right)} \right)}^{\frac{1}{\ga }}}\vgf(r)dr} 
  \EA \ee
  where $ c_{19} = \int_{0}^{\infty} {\pi ^{ - \frac{N}{2}}}{r^{\frac{N}{2}}}e^{-r}dr$ and $ \vgf(r) = c_{19}^{ - 1}{\pi ^{ - \frac{N}{2}}}{r^{\frac{N}{2}}}e^{ - r} $.
  Since $\ga \geq 1$ and $\int_{0}^{\infty}\vgf(r)dr=1$, thanks to Jensen's inequality we get 
$$ \BA{lll}
  \exp \left( {\delta {u^\ga }(x,t)} \right) &\leq \exp \left( {\delta {{\left( {\myint{0}{\infty}  {{c_{19}M_1}{{\left( { - \ln \left( { (4tr)^{\frac{1}{2}}  \wedge 2^{-1}} \right)} \right)}^{\frac{1}{\ga }}}\vgf(r)dr} } \right)}^\ga }} \right) \\ [3mm]
     & \leq \myint{0}{\infty}  {\exp \left( {\delta {{\left( {{c_{19}M_1}{{\left( { - \ln \left( { (4tr)^{\frac{1}{2}}  \wedge 2^{-1}} \right)} \right)}^{\frac{1}{\ga }}}} \right)}^\ga }} \right)\vgf(r)dr}. 
\EA $$
If $M_1=c_{19}^{-1}\gd^{-\frac{1}{\ga}}$ then
$$ \exp \left( {\delta {u^\ga }(x,t)} \right) \leq  \myint{0}{\infty}  {{\left( (4tr)^{\frac{1}{2}}  \wedge 2^{-1} \right)}^{-1}}\vgf(r)dr. $$
 Notice that $( (4tr)^{\frac{1}{2}}  \wedge 2^{-1})^{-1} \leq \left( {4tr} \right)^{ - \frac{1}{2}} + 2$, therefore
 \begin{equation*} \BA{ll}
 \exp \left( {\delta {u^\ga }(x,t)} \right) \le t^{ - \frac{1}{2}}\myint{0}{\infty} (4r)^{ - \frac{1}{2}}\vgf(r)dr  + 2,
\EA  \end{equation*}
which is \eqref{eu} with $c_{18}:=\int_{0}^{\infty} (4r)^{ - \frac{1}{2}}\vgf(r)dr$. \qed

\bth{nononlinear1} Assume $\gm \in \GTM^b_+(Q_T)$. There exists a positive constant $c_{20}=c_{20}(N)$ such that the unique solution $u$ of 
\bel{nononlinear1} \left\{ \BA{lll} \prt_t u - \Gd u &= \gm \qq &\text{in } Q_T \\
\phantom{----}
u &= 0 &\text{on } \prt \Gw \ti (0,T) \\
\phantom{----}
u &= 0 &\text{in }  \Gw 
\EA \right. \ee
satisfies $u(x,t) \leq c_{20}\BBW_{2d}[\gm](x,t)$ for every $(x,t) \in Q_T$.  
\es 
\Proof The unique solution of $(\ref{nononlinear1})$ is represented by (see \cite{QS})
\bel{repsol} u(x,t) = \myint{Q_t}{}G^\Gw(x,t-s,y)d\gm(y,s) \forevery (x,t) \in Q_T. \ee
Due to Fubini theorem, we obtain
$$ \BA{lll} u(x,t)& \leq \myint{Q_t}{}G(x-y,t-s)d\gm(y,s)\\
&=(4\pi)^{-\frac{N}{2}}\myint{Q_t}{}\left( \myfrac{N}{2}\myint{t-s}{\infty}\gt^{-\frac{N+2}{2}}d\gt \right)\left( \myint{\frac{|x-y|^2}{4(t-s)}}{\infty}e^{-r}dr  \right)d\gm(y,s) \\
&=2^{-N-1}N\pi^{-\frac{N}{2}} \myint{0}{\infty}\myint{0}{\infty}\gt^{-\frac{N+2}{2}} e^{-r}\gm(B_{\sqrt {4r\tau}}(x) \ti (t-\gt,t) )dr d\gt \\
&= I_1 + I_2
\EA $$
where 
$$  I_1:= 2^{-N-1}N\pi^{-\frac{N}{2}} \myint{0}{\infty}\myint{0}{\frac{1}{2}}\gt^{-\frac{N+2}{2}} e^{-r}\gm(B_{\sqrt {4r\tau}}(x) \ti (t-\gt,t) )dr d\gt, $$
$$ I_2:= 2^{-N-1}N\pi^{-\frac{N}{2}} \myint{0}{\infty}\myint{\frac{1}{2}}{\infty}\gt^{-\frac{N+2}{2}} e^{-r}\gm(B_{\sqrt {4r\tau}}(x) \ti (t-\gt,t) )dr d\gt. $$ 
By change of variables, we deduce
$$ \BA{lll} I_1 \leq 2^{-\frac{N}{2}}N\pi^{-\frac{N}{2}}(1-e^{-\frac{1}{2}})\BBW_{\infty}[\gm](x,t), \\[3mm]
 I_2  \leq  N\pi^{-\frac{N}{2}}\left( \myint{\frac{1}{2}}{\infty}r^{\frac{N}{2}}e^{-r}dr\right)\BBW_{\infty}[\gm](x,t).
\EA $$
Therefore 
\bel{f3} u(x,t) \leq c_{21} \BBW_{\infty}[\gm](x,t) \forevery (x,t) \in Q_T \ee
where $c_{21}=2^{-\frac{N}{2}}N\pi^{-\frac{N}{2}}(1-e^{-\frac{1}{2}})+N\pi^{-\frac{N}{2}}( \int_{\frac{1}{2}}^{\infty}r^{\frac{N}{2}}e^{-r}dr)$. By combining $(\ref{f3})$ and the estimate $ \BBW_{\infty}[\gm](x,t) < \frac{2^N}{2^N-1}\BBW_{2d}[\gm](x,t) $, we finish the proof. \qed

\bth{nononlinear} Assume $q \geq 1$, $\gd>0$, $\ga \geq q$, $\gb \in [\frac{q-1}{q},1)$,  $\gw \in \GTM_+^b(\Gw)$ and $\gm \in \GTM_+^b(Q_T)$. There exist  $M_1=M_1(N,\ga,\gd)$, $M_2=M_2(N,\gb,\gd)$ and $c_{22}=c_{22}(N,T,\Gw,d,\gd)$ such that if  $\norm{\BBM_{\ga,\infty}^{1}[\gw]}_{L^\infty(\BBR^N)} \leq M_1$ and $\norm{\BBM_{\gb,\infty}^{2}[\gm]}_{L^\infty(\BBR^{N+1})} \leq M_2$ then the unique solution $u$ of
\bel{nononlinear}\left\{ \BA{lll}
       u_t - {\Delta }u &= \mu \qq &\text{in } Q_T\\ 
      \phantom{-,,,,,}
       u  &= 0 &\text{on } \prt \Gw \ti (0,T) \\
      \phantom{-,} 
       u(.,0) &=\omega &\text{in }  \Gw  
\EA \right.\ee   
satisfies
\bel{est1} u \leq \BBG[\gw] + c_{20}\BBW_{2d}[\gm] \q \text{in } Q_T, \ee
\bel{est2} \myint{Q_T}{}\exp( \delta  u^q(x,t) )dxdt \leq c_{22} \ee
where  $c_{20}$ is the constant in \rth{nononlinear1}.  When $\ga=\frac{1}{1-\gb}=q$ then $c_{22}$ is independent of $\gd$.
\es
\Proof Let $v$ and $w$ be the solution of  \eqref{nononlinear0} and \eqref{nononlinear1} in $Q_T$ respectively. The  function $u:=v + w$ is the unique solution of \eqref{nononlinear} in $Q_T$. Hence estimate \eqref{est1}) follows from \rth{nononlinear0} and \rth{nononlinear1}. 

We next prove $(\ref{est2})$. By taking into account the fact that $e^{a+b} \leq 2^{-1}(e^{2a}+e^{2b})$ for every $a,b \in \BBR$, from \rth{nononlinear0} and \rth{nononlinear1}, we get 
\bel{L0} \BA{lll} \myint{Q_T}{}\exp\left( \delta  u^q \right)dxdt &\leq \myint{Q_T}{}\exp\left[ \delta  2^{q-1}(v^q+w^q) \right]dxdt \\[3mm]
&\leq \myfrac{1}{2}\myint{Q_T}{}\left(\exp\left[ \delta  (2v)^q \right] + \exp\left[ \delta  (2w)^q \right]\right)dxdt \\[3mm]
&\leq \myfrac{1}{2}\myint{Q_T}{}(\exp\left[ \delta  (2v)^q \right] + \exp\left[ \delta (2c_{20} \BBW_{2d}[\gm])^q \right])dxdt. 
\EA \ee
Next we set $M_1=2^{-1}c_{19}^{-1}\gd^{-\frac{1}{\ga}}$. It follows from \rth{nononlinear0} that if  $\norm{\BBM_{\ga,\infty}^1[\gw]}_{L^\infty(\BBR^N)} \leq M_1$ then 
\bel{L1} \myint{Q_T}{}(\exp\left( \delta  (2v)^\ga \right)dxdt < c_{23} \ee
where $c_{23}=c_{23}(N,T,\Gw)$.
Put $M_2=2^{-1}c_{20}^{-1}\gd_2^{1-\gb}\gd^{\gb-1}$ where $\gd_2$ is the constant in \rth{Hexp}. By \rth{Hexp}, if $\norm{\BBM_{\gb,\infty}^2[\gm]}_{L^\infty(\BBR^{N+1})}\leq M_2$ then
\bel{L2} \myint{Q_T}{}\exp\left[ \delta (2c_{20} \BBW_{2d}[\gm])^{\frac{1}{1-\gb}} \right]dxdt <c_{24} \ee
where $c_{24}=c_{24}(N,\gb,d)$.

Since $\ga \geq q$ and $\frac{1}{1-\gb} \geq q$, by combining Young inequality with \eqref{L0}, \eqref{L1} and \eqref{L2}, we derive \eqref{est2}. Notice that if $\ga=\frac{1}{1-\gb}=q$ then $c_{22}$ is independent of $\gd$. \qed
\section{Applications}
Let $\Gw$ be a bounded domain with $C^2$ boundary. This section is devoted to the proof of \rth{exist2} and \rth{exist3}.
\subsection{Equations with absorption terms}
 We first study the existence and uniqueness of solution to the following problem
\bel{absorp}\left\{ \BA{lll}
       u_t -  \Gd u + sign(u)g_\ell(u)&= \mu \qq &\text{in } Q_T \\ 
      \phantom{---------,,}
       u  &= 0 &\text{on } \prt \Gw \ti (0,T) \\ 
      \phantom{------,,,,,}
       u(.,0) &=\omega &\text{in }  \Gw  \\  
\EA \right.\ee   
where $g_\ell$ is define as in $(\ref{g})$ with $a>0$, $q \geq 1$, $\ell  \geq 1$. 

\bth{exist1} Assume $q\geq 1$, $a>0$, $\ga \geq q$, $\gb \in [\frac{q-1}{q},1)$,  $\gw \in \GTM^b(\Gw)$ and $\gm \in \GTM^b(Q_T)$. There exist positive constants $M_1=M_1(a,\ga,N)$ and $M_2=M_2(a,\gb,N)$ such that if $\norm{\BBM_{\ga,\infty}^{1}[\gw^{\pm}]}_{L^\infty(\BBR^N)} \leq M_1$ and $\norm{\BBM_{\gb,\infty}^2[\gm^{\pm}]}_{L^\infty(\BBR^{N+1})}\leq M_2$ then  the problem $(\ref{absorp})$ admits a unique solution $u$ satisfying $e^{a|u|^q} \in L^1(Q_T)$.
\es
\Proof {\it Step 1: Uniqueness.} If $u_1$ and $u_2$ are two solution of $(\ref{absorp})$ with the same data $(\gw,\gm) \in \GTM^b(\Gw) \ti \GTM^b(Q_T)$  then $u=u_1-u_2$ is a solution to problem
$$ \left\{ \BA{lll}
       u_t -  \Gd u + sign(u_1)g_\ell(u_1)-sign(u_2)g_\ell(u_2)&= 0 \q&\text{in } Q_T \\ 
      \phantom{--------------,,,,,,,,,,}
       u  &= 0 &\text{on } \prt \Gw \ti (0,T) \\ 
      \phantom{------------,,,,,,,,,,}
       u(.,0) &=0 &\text{in }  \Gw.  \\  
\EA \right.$$  
By \cite[Lemma 1.6 iii)]{MV0}, for every nonnegative function $\zeta \in X(Q_T)$,
\bel{uniq} \BA{lll} \myint{Q_T}{}-(\zeta_t +\Gd \zeta)|u|dxdt \\
\phantom{--}
 +  \myint{Q_T}{}\zeta \sign(u_1-u_2)(sign(u_1)g_\ell(u_1)-sign(u_2)g_\ell(u_2))dxdt \leq 0. \EA \ee 
Since the second term on the right-hand side in \eqref{uniq} is nonnegative, it follows that $\int_{Q_T}-(\zeta_t +\Gd \zeta)|u|dxdt \leq 0$. By choosing $\zeta=\psi$ which satisfies 
$$ \left\{ \BA{lll}
       -\psi_t -  \Gd \psi & =1 \qq &\text{in } Q_T \\ 
      \phantom{---,,}
       \psi  &= 0 &\text{on } \prt \Gw \ti (0,T) \\ 
      \phantom{--}
       \psi(.,T) &=0 &\text{in }  \Gw  \\  
\EA \right.$$  
we deduce that $u \equiv 0$, namely $u_1 \equiv u_2$.

It remains to deal with the question of existence. \medskip 

\noindent {\it Step 2: Approximating solutions}. Put $\gw_{1,n}=\rho_n * \gw^+$, $\gw_{2,n}=\rho_n * \gw^- $, $\gw_n=\gw_{1,n}-\gw_{2,n}$, $\gm_{1,n}=\tl \rho_n * \gm^+$, $\gm_{2,n}=\tl \rho_n * \gm^-$, $\gm_n=\gm_{1,n}-\gm_{2,n}$ where $\{\rho_n\}$ and $\{\tl \rho_n\}$ are sequences of mollifiers in $\BBR^N$ and $\BBR^{N+1}$ respectively. We may assume that $\gw_{i,n} \in C_c^\infty(\Gw)$ and $\gm_{i,n} \in C_c^\infty(Q_T)$ for every $n$ and $i=1,2$.  For each $n>0$, let $u_n$, $u_{i,n}$, $v_{i,n}$ ($i=1,2$) be respectively solutions to
\bel{absorp-n}\left\{ \BA{lll}
       (u_n)_t -  \Gd u_n + sign(u_n)g_\ell(u_n)&= \mu_n \qq &\text{in } Q_T \\ 
      \phantom{----------,,,,,,}
       u_n  &= 0 &\text{on } \prt \Gw \ti (0,T) \\ 
      \phantom{----------,}
       u_n(.,0) &=\gw_n &\text{in }  \Gw  \\  
\EA \right.\ee 

\bel{absorp-in}\left\{ \BA{lll}
       (u_{i,n})_t -  \Gd u_{i,n} + g_\ell(u_{i,n})&= \mu_{i,n} \qq &\text{in } Q_T \\ 
      \phantom{-------,,,,,,}
       u_{i,n}  &= 0 &\text{on } \prt \Gw \ti (0,T) \\ 
      \phantom{------,,,,}
       u_{i,n}(.,0) &=\gw_{i,n} &\text{in }  \Gw  \\  
\EA \right.\ee 

\bel{in}\left\{ \BA{lll}
       (v_{i,n})_t -  \Gd v_{i,n} &= \mu_{i,n} \qq &\text{in } Q_T \\ 
      \phantom{---,,,,,}
       v_{i,n}  &= 0 &\text{on } \prt \Gw \ti (0,T) \\ 
      \phantom{--,,,}
       v_{i,n}(.,0) &=\gw_{i,n} &\text{in }  \Gw  \\  
\EA \right.\ee 
By the maximum principle, $ -v_{2,n} \leq -u_{2,n} \leq u_n \leq u_{1,n} \leq v_{1,n} $
in $Q_T$. Therefore, $\abs{u_n} \leq \max\{u_{1,n},u_{2,n}\} \leq max\{v_{1,n},v_{2,n}\}$ in $Q_T$.
\medskip

\noindent{\it Step 3: End of proof.} Since $\{\gw_n\}$ and $\{\gm_n\}$ converge weakly to $\gw$ and $\gm$ respectively,  there exists a function $u$ and a subsequence, still denoted by $\{u_n\}$, such that $\{u_n\}$ and $\{g(u_n)\}$ converge to $u$ and $g(u)$ a.e. in $Q_T$. 

By \cite{MV1}, for any $p \in [1,\frac{N+2}{N})$, there exists a constant $c_{25}=c_{25}(\Gw,T,p)$ such that
$$ \norm{v_{i,n}}_{L^p(Q_T)} \leq c_{25}(\norm{\gm_{i,n}}_{L^1(Q_T)}+\norm{\gw_{i,n}}_{L^1(\Gw)}) \leq c_{25}(\norm{\gm_{i}}_{\CM(Q_T)}+\norm{\gw_{i}}_{\CM(\Gw)}),$$ 
from which it follows that
$$ \norm{u_n}_{L^p(Q_T)} \leq  c_{25}(\norm{\gm_{i}}_{\CM(Q_T)}+\norm{\gw_{i}}_{\CM(\Gw)}).$$ 
Therefore, due to Holder inequality, the sequence $\{u_n\}$ is equi-integrable. By Vitali theorem, the sequence $\{u_n\}$ converges to $u$ in $L^1(Q_T)$.
 
Notice that if $\norm{\BBM_{\ga,\infty}^{1}[\gw^{+}]}_{L^\infty(\BBR^N)} \leq M_1$ for some $M_1>0$ then for every $n \in \BBN$, $\norm{\BBM_{\ga,\infty}^{1}[\gw_{1,n}]}_{L^\infty(\BBR^N)} \leq M_1$. Indeed, for every $x \in \BBR^N$ and $s>0$, by Fubini theorem, we get
\bel{moll} \BA{lll}
 {\gw _{1,n}}({B_s}(x)) &=\myint{{B_s}(x)}{} {\myint{\mathbb{R}^N}{} {{\rho _n}(y - z)d\omega^+ (z)} dy} \\ [3mm]
    &=  \myint{\mathbb{R}^N}{} \myint{\mathbb{R}^N}{}{{\chi _{{B_s}(x - z)}}(y){\rho _n}(y)dyd\omega^+ (z)}   \\ [3mm]
    &= \myint{\mathbb{R}^N}{}\myint{\mathbb{R}^N}{} {{\chi _{{B_s}(x - z)}}(y)d\omega^+ (z){\rho _n}(y)dy}  \\[3mm]
    &= \myint{\mathbb{R}^N}{} {\omega^+ \left( {{B_s}(x - y)} \right){\rho _n}(y)dy}.
 \EA \ee
 Since ${\omega^+ \left( {{B_s}(x - y)} \right)}\le M_1s^Nh_{1,\ga}(s)$, we get $\omega_{1,n}(B_s(x))\le M_1s^Nh_{1,\ga}(s)$. Hence $\norm{\BBM_{\ga,\infty}^{1}[\gw_{1,n}]}_{L^\infty(\BBR^N)} \leq M_1$. Similarly, if $\norm{\BBM_{\gb,\infty}^2[\gm^{\pm}]}_{L^\infty(\BBR^{N+1})}\leq M_2$ for some $M_2>0$ then for every $n \in \BBN$, $\norm{\BBM_{\gb,\infty}^2[\gm_{i,n}]}_{L^\infty(\BBR^{N+1})}\leq M_2$.

Therefore, by setting $M_1=2^{-\frac{\ga+1}{\ga}}c_{19}^{-1}a^{-\frac{1}{\ga}}$ and $M_2:=2^{\gb-2}c_{20}^{-1}\gd_2^{1-\gb}a^{\gb-1}$, by \rth{nononlinear}, if $\norm{\BBM_{\ga,\infty}^{1}[\gw^{\pm}]}_{L^\infty(\BBR^N)} \leq M_1$ and $\norm{\BBM_{\gb}^2[\gm^{\pm}]}_{L^\infty(\BBR^{N+1})}\leq M_2$ then $ \int_{Q_T}\exp(2a v_{i,n}^q)dxdt \leq c_{26}$ where $c_{26}=c_{26}(N,T,\Gw,d,a)$. It follows that $\int_{Q_T}\exp(2a |u_n|^q)dxdt \leq c_{26}$. 
Consequently, $\{sign(u_n)g_\ell(u_n)\}$ is equi-integrable. Hence, by Vitali theorem, up to a subsequence, $\{sign(u_n)g_\ell(u_n)\}$ converges to $sign(u)g_\ell(u)$ in $L^1(Q_T)$. 

The solution $u_n$ satisfies, for every $\gz \in X(Q_T)$, 
\bel{solun} \myint{Q_T}{}(-u_n(\gz_t + \Gd \gz) + sign(u_n)g_\ell(u_n))\gz)dxdt = \myint{Q_T}{}\gz d\gm_n + \myint{\Gw}{}\gz(.,0)d\gw_n\ee
By letting $n \to \infty$ in $(\ref{solun})$, we deduce that $u$ is a solution to $(\ref{absorp})$. \qed

\blemma{apriori} Assume $\gw \in \GTM_+^b(\Gw)$ and $\gm \in \GTM_+^b(Q_T)$. Let $\{\gw_n\}$, $\{\gm_n\}$ and $\{u_n\}$ be defined as in step 1 of the proof of \rth{exist1}. There holds
\bel{normg} \norm{g_\ell(u_n)}_{L^1(Q_T)} \leq \norm{\gm}_{\GTM(Q_T)}+\norm{\gw}_{\GTM(\Gw)} \forevery n \in \BBN. \ee
\es
\Proof For any $k>0$, define $T_k(s)=\min\{k,\max\{-k,s\}\}$ for $s \in \BBR$ and $\ovl T_k(s)=\int_{0}^{s}T_k(\gs)d\gs$. For any $n \in \BBN$, $\ge>0$ , the function $\ge^{-1}T_\ge(u_{i,n})$ ($i=1,2$) can be employed as a test function for the problem $(\ref{absorp-in})$, i.e.
$$ \BA{lll} \myint{Q_T}{}\ge^{-1}(\ovl T_\ge(u_{i,n}))_tdxdt + \ge^{-1}\myint{Q_T}{}|\nabla T_\ge(u_{i,n})|^2dxdt \\
\phantom{------}
+ \myint{Q_T}{}g_\ell(u_{i,n})\ge^{-1}T_\ge(u_{i,n})dxdt = \myint{Q_T}{}\ge^{-1}T_\ge(u_{i,n})\gm_{i,n}dxdt. \EA $$
Since 
$$ \BA{lll}\myint{Q_T}{}(\ge^{-1}\ovl T_\ge(u_{i,n}))_tdxdt &=\myint{\Gw}{}\ge^{-1}\ovl T_\ge(u_{i,n}(T))dx - \myint{\Gw}{}\ge^{-1}\ovl T_\ge(\gw_{i,n})dx \\
&\geq - \norm{\gw_{i,n}}_{L^1(\Gw)}, \EA $$
it follows that
$$ \BA{lll}\myint{Q_T}{}g_\ell(u_{i,n})\ge^{-1}T_\ge(u_{i,n})dxdt &\leq \norm{\gm_{i,n}}_{L^1(Q_T)} + \norm{\gw_{i,n}}_{L^1(\Gw)}\\ 
&\leq \norm{\gm}_{\GTM(Q_T)}+\norm{\gw}_{\GTM(\Gw)}. \EA $$
By letting $\ge \to 0$, we derive $(\ref{normg})$. \qed \medskip

\noindent{\bf  Proof of \rth{exist2}.} For each $k>0$, $n \in \BBN$, denote by $u:=u^{f_1,f_2}_{k,n}$ the solution of
\bel{fn}\left\{ \BA{lll}
       u_t -  \Gd u + sign(u)g_\ell(u)&= \tl \rho_n * \mu + \tl \rho_n *( T_k(f_2)) &\text{in } Q_T \\ 
      \phantom{---------,}
       u  &= 0 &\text{on } \prt \Gw \ti (0,T) \\ 
      \phantom{-------,,}
       u(.,0) &=\rho_n * \gw + \rho_n * (T_k(f_1))&\text{in }  \Gw  \\  
\EA \right.\ee   
where $\{\rho_n\}$ and $\{ \tl \rho_n\}$ are sequences of mollifiers in $\BBR^N$ and $\BBR^{N+1}$ respectively.  Let $u:=u^{f_1,f_2}_{\pm,k,n}$ is the solution of
\bel{fn+-}\left\{ \BA{lll}
        u_t -  \Gd u + g_\ell(u)&= \tl \rho_n *( \mu^\pm) + \tl \rho_n *( T_k(f_2^\pm))\q &\text{in } Q_T \\ 
      \phantom{------,,}
       u  &= 0 &\text{on } \prt \Gw \ti (0,T) \\ 
      \phantom{----,,}
       u(.,0) &=\rho_n *( \gw^\pm) + \rho_n * (T_k(f_1^\pm))&\text{in }  \Gw.  \\
\EA \right.\ee   
By the comparison principle, $-u^{f_1,f_2}_{-,k,n} \leq u^{f_1,f_2}_{k,n} \leq u^{f_1,f_2}_{+,k,n}$ for any $n$. Using similar argument as in \rth{exist1}, we deduce that there exist  $M_1=M_1(N,\ga,a)$ and $M_2=M_2(N,\gb,a)$ such that if $\norm{\BBM_{\ga,\infty}^{1}[\gw^{\pm}]}_{L^\infty(\BBR^N)} < M_1$ and $\norm{\BBM_{\gb,\infty}^{2}[\gm^{\pm}]}_{L^\infty(\BBR^{N+1})}<M_2$, there holds
$$ \myint{Q_T}{}\exp(2a |u^{f_1,f_2}_{\pm,k,n}|^q)dxdt \leq c_{27} $$
where $c_{27}=c_{27}(N,T,\Gw,\gb,a,d,k)$. Hence we can find a subsequence, still denoted by  $\{u^{f_1,f_2}_{\pm,k,n}\}$, and a function $u^{f_1,f_2}_{\pm,k}$, such that $\{u^{f_1,f_2}_{\pm,k,n}\}$ and $\{ g_\ell(u^{f_1,f_2}_{\pm,k,n})\}$ converge to $u^{f_1,f_2}_{\pm,k}$ and $g_\ell(u^{f_1,f_2}_{\pm,k})$ respectively in $L^1(Q_T)$ as $n \to \infty$. Therefore, $u^{f_1,f_2}_{\pm,k}$ is the solution of $(\ref{absorp})$ with $\gm$ replaced by $\gm^{\pm} + T_k(f_2^{\pm})$  and $\gw$ replaced by $\gw^{\pm} + T_k(f_1^{\pm})$. By a similar argument, we can show that there exists a unique solution $u^{f_1,f_2}_{k}$  of $(\ref{absorp})$ with $\gm$ replaced by $\gm + T_k(f_2)$ and $\gm$ replaced by $\gw + T_k(f_1)$. Moreover, by the comparison principle, $-u^{f_1,f_2}_{-,k} \leq u^{f_1,f_2}_{k} \leq u^{f_1,f_2}_{+,k}$ and the sequences $\{u^{f_1,f_2}_{\pm,k}\}$ are increasing with respect to $k$. Denote $u^{f_1,f_2}_{\pm}:=\lim_{k \to \infty}u^{f_1,f_2}_{\pm,k}$.Thanks to \rlemma{apriori}that for every $k>0$,
$$ \myint{Q_T}{}g_\ell (u^{f_1,f_2}_{\pm,k})dxdt \leq \norm{\gw}_{\CM(\Gw)}+ \norm{\gm}_{\CM(Q_T)}+\norm{f_1}_{L^1(\Gw)} + \norm{f_2}_{L^1(Q_T)}. $$
Therefore, by monotone convergence theorem, $\{g_\ell(u^{f_1,f_2}_{\pm,k})\}$ converges to $g_\ell (u^{f_1,f_2}_{\pm})$ in $L^1(Q_T)$.

Since $-u^{f_1,f_2}_{-,k} \leq u^{f_1,f_2}_{k} \leq u^{f_1,f_2}_{+,k}$, it follows that $g_\ell(u^{f_1,f_2}_k) \leq g_\ell(u^{f_1,f_2}_{-,k}) + g_\ell(u^{f_1,f_2}_{+,k})$. Therefore the sequence $\{ \gm + T_k(f_2) - sign(u^{f_1,f_2}_k)g_\ell(u^{f_1,f_2}_k) \}$ is bounded in $\GTM(Q_T)$. Notice that the sequence $\{ \gw + T_k
(f_1)\}$ is also bounded in $\GTM(\Gw)$. Hence, up to a subsequence, $\{ u^{f_1,f_2}_k\}$ converges to a function $u^{f_1,f_2}$ in $L^1(Q_T)$ and a.e. in $Q_T$.  Moreover, by dominated convergence theorem, the sequence $\{sign(u^{f_1,f_2}_k)g_\ell(u^{f_1,f_2}_{k})\}$ converges to $sign(u^{f_1,f_2})g_\ell(u^{f_1,f_2})$ in $L^1(Q_T)$ as $k \to \infty$. By passing to the limite, we deduce that $u^{f_1,f_2}$ is a solution of $(\ref{absorp-f})$. The uniqueness is obtained by using similar argument as in \rth{exist1}. \qed
\subsection{Equations with source terms}
In this section we deal with the existence of solutions to problem \eqref{source}. \medskip

\noindent{\bf Proof of \rth{exist3}}. Let $u_0$ be a solution of
$$ \left\{ \BA{lll} \prt_t u_0 - \Gd u_0 &= \gm \qq &\text{in } Q_T\\
\phantom{----,}
u_0 &= 0 &\text{on } \prt \Gw \ti (0,T) \\
\phantom{--,}
u_0(.,0) &= \gw &\text{in } \Gw.  \\
\EA \right.$$
For each $n \in \BBN$, let $u_{n+1}$ be a solution of 
$$ \left\{ \BA{lll} \prt_t u_{n+1} - \Gd u_{n+1} &= g_\ell(u_n)+ \gm \qq &\text{in } Q_T\\
\phantom{-----}
u_{n+1} &= 0 &\text{on } \prt \Gw \ti (0,T) \\
\phantom{---}
u_{n+1}(.,0) &= \gw &\text{in } \Gw  \\
\EA \right.$$
namely, for every $\zeta \in X(Q_T)$,
\bel{soln} -\myint{Q_T}{}u_{n+1}(\gz_t + \Gd \gz)dxdt = \myint{Q_T}{}g_\ell(u_n)\zeta dx dt + \myint{Q_T}{}\gz d\gm + \myint{\Gw}{}\gz(.,0)d\gw\ee

We need the following lemma
\blemma{lem} Under the assumptions of \rth{exist3}, there exist positive constants $b_0=b_0(N,d,\ell,q) \in (0,1]$, $M_1=M_1(N,a,\ga,q,\ell,d)$, $M_2=M_2(N,a,\gb,q,\ell,d)$ such that if $\norm{\BBM_{\ga,\infty}^1[\gw]}_{L^\infty(\BBR^N)} \leq M_1$ and $\norm{\BBM_{\gb,\infty}^2[\gm]}_{L^\infty(\BBR^{N+1})} \leq M_2$ then
\bel{estm} u_n \leq \BBG[\gw] + c_{20} \BBW_{2d}[\gm] + c_{20}b_0 \forevery n \in \BBN. \ee
\bel{exploc} g_\ell(2 \BBG[\gw]+2c_{20} \BBW_{2d}[\gm]+2c_{20}b_0) \in L^1(Q_T) \ee
where $c_{20}$ is the constant in \rth{nononlinear1}.
\es
\noindent{\it Proof of \rlemma{lem}.} We prove $(\ref{estm})$ by recurrence. Indeed, $(\ref{estm})$ holds true if $n=0$ by the previous results. Assume now $(\ref{estm})$ holds true when $n=m$. We shall show that $(\ref{estm})$ remains true when $n=m+1$. By \rth{nononlinear}, 
\bel{goal1} u_{m+1} \leq \BBG[\gw] + c_{20} \BBW_{2d}[g_\ell(u_m)+\gm]=\BBG[\gw] + c_{20} \BBW_{2d}[\gm] + c_{20} \BBW_{2d}[g_\ell(u_m)]. \ee
Therefore, it's sufficient to prove that
\bel{goal2} \BBW_{2d}[g_\ell(u_m)] \leq b_0. \ee
Since $(\ref{estm})$ is valid when $n=m$, it follows that
$$ g_\ell(u_m) \leq 3^{-1}g_\ell(3 \BBG[\gw]) + 3^{-1}g_\ell(3c_{20} \BBW_{2d}[\gm]) + 3^{-1}g_\ell(3c_{20}b_0). $$
Keeping in mind that $g_\ell(s) \leq \vge^{\ell q} g_\ell(\vge^{-1}s)$ for every $s \geq 0$ and $\vge \in (0,1]$, we get for $\vge_1, \vge_2 \in (0,1]$, 
$$ g_\ell(u_m) \leq 3^{-1}\vge_1 g_\ell(3\vge_1^{-1} \BBG[\gw]) + 3^{-1}\vge_2g_\ell(3\vge_2^{-1}c_{20} \BBW_{2d}[\gm]) + 3^{-1}b_0^{\ell q}g_\ell(3c_{20}). $$
Hence
\bel{compt} \BA{lll} \BBW_{2d}[g_\ell(u_m)] \leq  3^{-1}\vge_1 \BBW_{2d}[g_\ell(3\vge_1^{-1} \BBG[\gw])] \\
\phantom{----}
+ 3^{-1}\vge_2 \BBW_{2d}[g_\ell(3\vge_2^{-1}c_{20} \BBW_{2d}[\gm])] + 3^{-1}b_0^{\ell q} g_\ell(3c_{20})\BBW_{2d}[1]. \EA \ee
We choose $b_0$ such that
\bel{b0} b_0^{\ell q} g_\ell(3c_{20})\BBW_{2d}[1] = b_0 \Llra b_0 = \left( \frac{16}{3} \gw_N d^3 g_\ell(3c_{20})\right)^{\frac{1}{\ell q-1}}\ee
where $\gw_N$ is the volume of unit ball in $\BBR^N$. \medskip

\noindent{\it Step 1: We show that if $\vge_1$ is small enough then}
\bel{choice1} \vge_1 \BBW_{2d}[g_\ell(3\vge_1^{-1} \BBG[\gw])] \leq b_0.\ee
When $q=\ga$, by \rth{nononlinear0}, if $\norm{\BBM_\ga^{1}[\gw]}_{L^\infty(\BBR^N)} \leq 3^{-1}c_{19}^{-1}a^{-\frac{1}{\ga}}\vge_1$, then by \rth{nononlinear0},  for any $(y,\gt) \in Q_T$,
\bel{es} \BA{ll} g_\ell(3\vge_1^{-1} \BBG[\gw](y,\gt)) \leq \exp(a3^q\vge_1^{-q} (\BBG[\gw](y,\gt))^q) 
\leq c_{18}\gt^{-\frac{1}{2}} + 2 \EA \ee
where $c_{18}$ is the constant in \rth{nononlinear0}. Therefore, for and $s \geq 0$ and fixed $(x,t) \in Q_T$,
\bel{locin} \BA{lll} \myint{Q_s(x,t)}{}g_\ell(3\vge_1^{-1} \BBG[\gw](y,\gt))dyd\gt &= \myint{B_s(x)}{}\myint{(t-\frac{s^2}{2}) \lor 0}{t+\frac{s^2}{2}}g_\ell(3\vge_1^{-1} \BBG[\gw](y,\gt))dyd\gt \\ [5mm]
&\leq 2c_{18}\gw_Ns^{N+1} + \gw_Ns^{N+2}.
\EA \ee
Consequently, 
\bel{es1} \BA{ll} \BBW_{2d}[g_\ell(3\vge_1^{-1} \BBG[\gw])](x,t) \leq 4c_{18}\gw_N d^2+4\gw_Nd^3=:c_{28}
\EA \ee
Thus, if $ \vge_1c_{28} \leq b_0$, namely $\vge_1 \leq  c_{28}^{-1}b_0$, then $(\ref{choice1})$ holds true. When $q<\ga$, by Young inequality and \rth{nononlinear0}, if $\norm{\BBM_\ga^{1}[\gw]}_{L^\infty(\BBR^N)} \leq 3^{-\frac{q}{\ga}}c_{19}^{-1}a^{-\frac{1}{\ga}}\vge_1$ then $\BBW_{2d}[g(3\vge_1^{-1} \BBG[\gw])](x,t) \leq c_{28}e^{a3^q}$. Hence if $ \vge_1 c_{28}e^{a3^q} \leq b_0$, namely $\vge_1 \leq c_{28}^{-1}e^{-a3^q}b_0 $
then $(\ref{choice1})$ holds true. Thus, by putting $\vge_1=(c_{28}^{-1}e^{-a3^q}b_0) \wedge 1$ and $M_1=3^{-1}c_{19}^{-1}a^{-\frac{1}{\ga}}((c_{28}^{-1}e^{-a3^q}b_0) \wedge 1)$, we derive $(\ref{choice1})$ for every $\ga \geq q$. \medskip

\noindent{\it Step 2:  We show that if $\vge_2$ small enough then }
\bel{choice3} \vge_2 \BBW_{2d}[g_\ell(3\vge_2^{-1}c_{20} \BBW_{2d}[\gm])] \leq b_0. \ee
When $q=(1-\gb)^{-1}$, thanks to \rth{Hexp}  if 
$$ a3^q\vge_2^{-q}c_{20}^q \leq \gd_2 M_2^{-\frac{1}{1-\gb}} \Llra M_2 \leq a^{\gb-1}3^{-1}c_{20}^{-1}\vge_2\gd_2^{1-\gb}$$ 
then for any $0<s <2d$, there holds
\bel{loc} \myint{Q_s(x,t)}{}g_\ell(3\vge_2^{-1}c_{20} \BBW_{2d}[\gm]) dyd\tau \leq c_{29}. \ee
and 
$$ \BA{lll} \BBW_{2d}[g_\ell(3\vge_2^{-1}c_{20} \BBW_{2d}[\gm])] &< \BBW_{2d}[\exp(a3^q\vge_2^{-q}c_{20}^q (\BBW_{2d}[\gm])^q)] \\[3mm]
&\leq \BBW_{2d}[\exp(\gd_2 M_2^{-q}(\BBW_{2d}[\gm])^q)] \leq c_{30} \EA $$
where $c_{i}=c_{i}(N,\gb,d)$ with $i=29,30$. Hence it's sufficient to choose $\vge_2$ such that $\vge_2 c_{30} \leq b_0$, i.e. $\vge_2 \leq c_{30}^{-1}b_0$.
When $q<(1-\gb)^{-1}$, by Young inequality and \rth{Hexp}, if
$$ a3^q c_{20}^q\vge_2^{-\frac{1}{1-\gb}}\leq \gd_2 M_2^{-\frac{1}{1-\gb}} \Llra M_2 \leq a^{\gb-1}3^{-q(1-\gb)}c_{20}^{-q(1-\gb)}\vge_2\gd_2^{1-\gb} $$
then $\BBW_{2d}[g_\ell(3\vge_2^{-1}c_{20} \BBW_{2d}[\gm])] \leq c_{30}e^{a3^q c_{20}^q} $. Therefore, if 
$$ \vge_2 c_{30} e^{a 3^qc_{20}^q} \leq b_0 \Llra \vge_2 \leq c_{30}^{-1}e^{-a 3^q c_{20}^{-q}}b_0 $$
then $(\ref{choice3})$ follows. Thus, by setting $\vge_2 = (c_{30}^{-1}e^{-a 3^q c_{20}^{-q}}b_0) \wedge 1$ and
$M_2=a^{\gb-1}3^{-1}c_{20}^{-1}\gd_2^{1-\gb}((c_{30}^{-1}e^{-a 3^q c_{20}^{-q}}b_0) \wedge 1)$, we obtain $(\ref{choice3})$. \medskip

\noindent{\it Step 3: End of proof.} By combining $(\ref{b0})$, $(\ref{choice1})$ and $(\ref{choice3})$, we deduce that if $\norm{\BBM_\ga^1[\gw]}_{L^\infty(\BBR^N)} \leq M_1$ and $\norm{\BBM_\gb^2[\gm]}_{L^\infty(\BBR^{N+1})} \leq M_2$ then $(\ref{goal2})$ and $(\ref{estm})$ hold true. 

Moreover, by convexity, for any $\gg \in (0,1)$, we have
$$ \BA{lll} g_\ell(2 \BBG[\gw]+2c_{20} \BBW_{2d}[\gm]+2c_{20}b_0)  &\leq \myfrac{\gg}{4(1+\gg)}g_\ell\left(\myfrac{8(1+\gg)}{\gg} \BBG[\gw]\right) \\
&+ \myfrac{\gg}{4(1+\gg)}g_\ell\left(\myfrac{8(1+\gg)}{\gg}c_{20} \BBW_{2d}[\gm]\right) \\
&+\myfrac{2+\gg}{2(1+\gg)}g_\ell\left(\myfrac{4(1+\gg)}{2+\gg}c_{20}b_0\right). \EA $$
We choose $\gg$ such that
$$ \myfrac{8(1+\gg)}{\gg} = 3(\vge_1^{-1} \wedge \vge_2^{-1}) \Llra \gg = \frac{8}{3(\vge_1^{-1} \wedge \vge_2^{-1})-8}. $$
Then 
$$ \BA{lll} g_\ell(2 \BBG[\gw]+2c_{20} \BBW_{2d}[\gm]+2c_{20}b_0)  \\
\phantom{--------}
\leq g_\ell(3\vge_1^{-1} \BBG[\gw]) + g_\ell(3\vge_2^{-1}c_{20} \BBW_{2d}[\gm]) + g_\ell(4c_{20}b_0), \EA $$
which, together with $(\ref{locin})$ and $(\ref{loc})$, implies $(\ref{exploc})$. \qed \medskip

{\it Let us now return to the proof of \rth{exist3}}. \medskip

By comparison principle, $\{u_n\}$ is increasing and converges to a function $u$ a.e. in $\Gw$.  Moreover, it follows from $(\ref{exploc})$ that the sequences $\{u_{n}\}$  and $\{g_\ell(u_{n})\}$  are uniformly bounded in $L^1(Q_T)$. Thanks to monotone convergence theorem, $\{u_{n}\}$  and $\{g_\ell(u_{n})\}$ converge to $u$ and $g_\ell(u)$ respectively in $L^1(Q_T)$. By letting $n \to \infty$ in $(\ref{soln})$, we derive that $u$ is a solution of $(\ref{source})$. \qed

\end{document}